\documentclass[review]{elsarticle}
\usepackage[hyphens,spaces,obeyspaces]{url}

\usepackage[a4paper, margin=2.5cm]{geometry}

\usepackage{lineno,hyperref}
\usepackage{hyperref}
\usepackage{xcolor}
\usepackage{amsmath}
\usepackage{amsfonts}
\usepackage{multirow}

\journal{IISE Transactions}

\begin{document}

\begin{frontmatter}

\title{Optimal Heterogeneous Asset Location Modeling for Expected Spatiotemporal Search and Rescue Demands using Historic Event Data}

\author{Zachary T. Hornberger}
\author{Bruce A. Cox}
\ead{bruce.cox@afit.edu}
\author{Brian J. Lunday}
\ead{brian.lunday@afit.edu}
\address{Department of Operational Sciences, Air Force Institute of Technology,\\ 2950 Hobson Way, Wright-Patterson AFB, OH 45433}
\fntext[myfootnote]{Corresponding author at \href{mailto:bruce.cox@afit.edu}{bruce.cox@afit.edu}}




\begin{abstract}
The United States Coast Guard is charged with the coordination of all search and rescue missions in maritime regions within the United States' purview.  Given the size of the Pacific Ocean and the limited resources available to respond to search and rescue missions in this region, the service seeks to posture its aligned fleet of maritime and aeronautical assets to reduce the expected response time for such missions.  Leveraging historic event records for the region of interest, we propose and demonstrate a two-stage solution approach.  In the first stage, we develop and apply a stochastic zonal distribution model to evaluate spatiotemporal trends for emergency event rates and corresponding response strategies to inform the probabilistic modeling of future rescue events respective locations, frequencies, and demands for support.  In the second stage, the results from the aforementioned analysis enable the parameterization and solution of a integer linear programming formulation to identify the best locations at which to station limited heterogeneous search and rescue assets.  Considering both the 50th and 75th percentile levels of forecast event and asset demand distributions using 7.5 years of historical event data, our models identify asset location strategies that respectively yield a 9.6\% and 17.6\% increase in coverage over current asset basing when allowing locations among current homeports and airports, as well as respective 67.3\% and 57.4\% increases in coverage when considering a larger set of feasible basing locations.
\end{abstract}


\end{frontmatter}


\section{Introduction}
The United States Coast Guard (USCG) is charged with coordinating search and rescue (SAR) missions in maritime regions within the United States' purview \cite{NatSARPlan2018}.  As the federal SAR coordinator for these emergencies, the USCG responds to SAR events with its own maritime and aeronautical assets, and coordinates as necessary with private, commercial, and Department of Defense assets in the vicinity of SAR events to provide supplementary support.  The time required for the USCG to receive an emergency notification, dispatch an asset, locate the individual(s) in distress, and provide aid can be the difference between life and death, so it is important to ensure response assets are well positioned to service expected demands.  

Organizationally, the command and control of USCG operations is divided among administrative districts, and each district is assigned to direct operations for a geographic area of responsibility (AOR).  Headquartered in Honolulu on the island of O'ahu, USCG District 14 is responsible for missions across most of the Pacific Ocean, a region over 12 million square nautical miles in size which the district supports with assets stationed among the Hawaiian Islands and Guam.  Despite having the largest geographic AOR, District 14 is allocated relatively few resources, a result of both resource priorities across USCG districts, and a historically lean budget authorization for a service responding to very high levels of demand \cite{CGBudgetArchive,CGBudget}.  Moreover, every USCG district must service multiple missions (e.g., search and rescue, drug interdiction, migrant interdiction, marine safety); although every search and rescue mission relates to saving lives and is of paramount importance, the assets within District 14 are needed to support a myriad of activities.

Given its limited available resources, District 14 seeks to posture its fleet strategically to minimize the expected response time to respond to emerging SAR missions.  The strategic posturing of the District 14 fleet is complicated by the size of the AOR, the need to consider a heterogeneous fleet of maritime and aeronautical assets, and the inherent uncertainty of both the emergency event locations and frequencies for SAR events.  Thus, it is first necessary to predict expected SAR demands and subsequently optimize both the locations of assets and allocation of expected demands among such assets within the District 14 fleet.

The work herein is informed by two subdisciplines from the technical literature:  spatiotemporal forecasting and location theory.

The two major threads of published spatiotemporal forecasting research relate to predicting either natural or man-made incidents.  The potential to leverage wind power via spatiotemporal forecasting is examined by Xie et al. \cite{xie2014short} and Tastu et al. \cite{tastu2011spatio}, whereas Tascikaraoglu et al. \cite{tascikaraoglu2016compressive} sought to predict meteorological qualities suitable to generate solar power.  In contrast, Madadgar and Moradkhani \cite{madadgar2014spatio} sought to predict droughts with Bayesian networks to mitigate agricultural impacts.  With a more somber view towards hazards, research by Hashimoto et al. \cite{hashimoto2013predicted} developed spatiotemporal predictions of nucleotide deposits in the forests in the vicinity of the Fukushima nuclear accident in Japan.  By comparison, most work related to man-made incidents pertains to predicting traffic congestion, whether vehicular \cite{frihida2002spatio,sirvio2008spatio} or pedestrian \cite{zhang2017deep}.  Other research seeks to predict the times and locations of deliberate, malicious behavior by humans, such as wildlife poaching \cite{gholami2017taking}, internet system attacks \cite{soldo2011blacklisting}, and crime \cite{wang2011spatio,wang2012spatio,wang2012spatiob}. Benigni and Furrer \cite{benigni2012spatio} sought to develop effective spatiotemporal predictions for insurgents using improvised explosive devices on a subset of roads in Baghdad, Iraq.  More closely related to the current study is work by Prasannakumar et al. \cite{prasannakumar2011spatio} to predict the times and locations of road accidents, albeit for incidents relatively restricted in their geography.  Marchione and Johnson \cite{marchione2013spatial} sought to predict maritime piracy incidents, but their consideration of events resulting from deliberate human behavior renders it ill-suited for the problem considered herein.  Instead, we seek to predict accidents (i.e., SAR events) resulting from  maritime human activity across an expansive region, and which may or not be related to natural phenomena such as rough seas and inclement weather.

Setting aside the question of weather as a strictly causal factor for SAR events, there do exist studies that forecast SAR demands that relate to the movement of people by maritime assets across an expansive region.  Most such works address the uncertainty of SAR demands by decomposing a maritime region of study into a rectangular mesh of zones and subsequently modeling the temporal nature of the aggregated, zonal demand volumes via a simulation-based approach using either Poisson or uniform distributions, with the spatial identification of SAR events occurring at the center of each zone \cite{afshartous2009us, akbari2017modular, karatas2017ilp}.  Related studies consider a more customized decomposition of a region using the zonal distribution model set forth by Azofra et al. \cite{azofra2007optimum}.  Ai et al. \cite{ai2015optimization} and Razi and Karatas \cite{razi2016multi} implemented variations of the zonal distribution model using static average historical demands as a deterministic demand volume for each zone, while considering zonal demands to occur at fixed points (i.e., the centroids) within the respective zones.  Given the size of the District 14 AOR, the disparity of demand volume across its geographic expanse, and the uncertainty inherent in SAR demand, we develop and implement a \textit{stochastic zonal distribution model} that leverages both the stochastic components of simulation-based methods and the customized zoning approach used by the aforementioned deterministic models.

Given the need to locate the District 14 fleet of maritime and aeronautical assets to service expected SAR demands, the second stage of our research is informed by the field of location theory.  Early works in location theory focused on where a decision maker should locate facilities \cite{toregas1971location,church1974maximal,daskin1981hierarchical} and how many facilities to emplace at each possible location \cite{hall1972application,berlin1974mathematical,baker1989non}.  Among these works are models fundamental to the discipline, including the set covering location problem (SCLP), maximal covering location problem (MCLP),  $p$-center problem, and $p$-median problem.  Given a limit on the range of a facility to cover demands, the SCLP identifies the minimal the number (or cost) of facilities required to cover all demands, whereas the MCLP maximizes the number of demands covered by a bounded number (or cost) of facilities.  Setting aside the restriction of a fixed covering range (radius in Cartesian space) and adopting a location-allocation framework, the $p$-center problem locates $p$ facilities among a set of possible sites and allocates a set of demands to the located facilities, seeking to minimize the maximum facility-to-assigned-demand distance.  Alternatively, the $p$-median problem likewise locates facilities and allocates demands to them, but it minimizes the total (or equivalently, the average) facility-to-assigned-demand distance.  Extensions to these early models consider facilities having capacities \cite{berman2007locating}, deterministic and stochastic demands \cite{wang2002algorithms,snyder2007stochastic}, multiple covering of demands \cite{daskin1988integration,batta1990covering}, hierarchical facility structures \cite{teixeira2008hierarchical}, and other variants \cite{snyder2006facility,mclay2009maximum}; an interested reader is referred to recent reviews by Jia et al. \cite{jia2007modeling} and Farahani et al. \cite{farahani2012covering}, or to any of four excellent books by Drezner and Hamacher \cite{drezner2001facility}, Daskin \cite{daskin2011network}, Laporte et al. \cite{laporte2015location}, or Church and Murray \cite{church2018location}, each of which provides a thorough treatment of the subject.  Most closely related to our problem is the $p$-median location problem, which we adapt for the heterogeneous District 14 fleet of assets, as well as the multiple objectives entailed by locating assets effectively yet minimizing changes to the existing enterprise.

Although it does not directly inform our work, it is appropriate to mention the related work in locating emergency response assets and, in a more limited stream of literature, locating assets that respond specifically to maritime emergencies.   A rich body of literature also exists for the location of emergency management system (EMS) assets.  Each of the recent surveys by Aringhieri et al. \cite{aringhieri2017emergency}, Gholami-Zanjani et al. \cite{gholami2018or}, and B{\'e}langer et al. \cite{belanger2019recent} is quite thorough.  Although EMS missions are conceptually related to the Coast Guard's SAR mission, both the distance metric (i.e., road distance) and the oft heterogeneous assets located for EMS modeling make it sufficiently different to preclude the adoption of existing models from the ground-based EMS literature to address the District 14 problem.  Previous studies regarding the placement of resources to respond to \textit{maritime} emergencies are more limited, and they typically restrict their examination to a single asset type or the operational range of a station.  Wagner and Radovilsky \cite{wagner2012optimizing} and Razi and Karatas \cite{razi2016multi} both developed allocation tools limited to the assignment of USCG boats, and Akbari et a. \cite{akbari2017modular} and Azofra et al. \cite{azofra2007optimum} developed models to locate maritime vessels.  In contrast, Karatas et al. \cite{karatas2017ilp} scoped their allocation plan to Turkish Coast Guard helicopters, whereas Afshartous et al. \cite{afshartous2009us} constructed a statistical model to locate air stations and concluded by suggesting their methodology could be applied to consider the allocation of resources at those stations.  Although these studies tend to include detailed operational requirements for their specific asset types within the models, they lack the general applicability to the location of heterogeneous assets necessary to address the problem examined herein.

This research makes two contributions to the field of resource location-and-allocation problems.  First, it introduces the stochastic zonal distribution model, an extension to the works of Azofra et al.  \cite{azofra2007optimum} and Razi and Karata \cite{razi2016multi} to develop spatiotemporal maritime incident forecasts using historic demand data.  Second, it formulates and applies an integer linear programming formulation that leverages the output of the stochastic zonal distribution analysis to optimize the placement of a fleet of heterogeneous assets such that the expected response time for future SAR events is minimized.  

The remainder of this paper is organized as follows.  Section 2 discusses the respective stochastic zonal distribution and the asset location-and-allocation models and their implementations.  Section 3 presents and discusses the results of applying these models to District 14's problem using seven years of recent historic SAR event data, and Section 4 summarizes major insights and provides recommendations for future research.

\section{Models and Solution Methodologies}
In this section, we present in Section \ref{sec:ModelParam} our methodology for parameterizing the model using historic demand data, after which we define in Section \ref{sec:LocModel} both the notation and math programming formulation to address USCG District 14's heterogeneous asset location problem. 

\subsection{Model Parameterization for Stochastic Demand} \label{sec:ModelParam}
In SAR operations, future emergencies are inherently uncertain in time and location.  This uncertainty presents a challenge to applying a prescriptive optimization model because it begets stochastic demands having decimal values, whereas the different assets to be located are, of course, integer-valued.  Moreover, the challenge of determining uncertain demand levels is amplified by the size of the AOR, the variability in demand across the AOR, and the variability of demand for different asset types.  To address these challenges, we extend the work of Razi and Karata \cite{razi2016multi} by developing the stochastic zonal distribution model.  This model addresses the spatiotemporal forecasting problem in three steps:  construction of SAR demand zones, stochastic characterization of zonal demand frequencies, and identification of the levels (i.e., type and numbers of each asset) of such SAR demands.

\subsection*{Step 1.  Hierarchical Clustering of Historical SAR Data}
First, the region of study is decomposed into manageable zones using a hierarchical $k$-means clustering algorithm.  This algorithm first sorts historical SAR data into mutually exclusive groups based upon the units which coordinated the emergency response and the assets operationally capable of responding.  For this study, SAR events were categorically classified as either \textit{boat/helicopter events} or \textit{cutter/airplane events}, where the former classification entails an event to which shorter range assets can adequately respond.  The islands with USCG boat stations are the Hawaiian islands of Kaua'i, O'ahu, and Maui as well as the island of Guam.  Due to the operational constraints of USCG boats and the tendency to combine boats and helicopters together in operations, SAR events within 50 nautical miles of these islands were classified as boat/helicopter events whereas emergencies beyond these boundaries were classified as cutter/airplane events.  Each group is then geographically clustered using a $k$-means technique.  The implementation of a $k$-means algorithm for determining zones within a region of study was inspired by the work of Razi and Karata \cite{razi2016multi}.

\subsection*{Step 2.  Probabilistic Representation of SAR Demand for each Zone}
Second, probabilistic representations of monthly SAR demands are constructed for each zone.  Given the stochastic nature of SAR emergencies, we identify probability distributions that accurately represent the emergence of SAR events over time.  Afshartous et al. \cite{afshartous2009us} and Akbari et al. \cite{akbari2017modular} identified SAR events as likely being Poisson-distributed within their simulation-based approach to locate SAR assets.  We hypothesize the historical SAR data for this study is likewise represented appropriately via Poisson processes, although the assumptions of independence and stationarity are evaluated.

\subsection*{Step 3.  Probabilistic Representation of SAR Event Asset Response Requirements}
Third, based upon the probabilistic representation of SAR demand for each zone, the corresponding response is modeled.  When the Coast Guard receives notification of an emerging crisis, a SAR team uses a system known as SAROPS to determine the appropriate assets to respond, given the weather, asset availability, and case details.  While this level of fidelity is beyond the scope of our research, we note the trends in rescue operations based on the geography of the region.  Reviewing event records, we classify response strategies in each zone and determine their respective frequencies of occurrence.  We bounded the response requirement to be no more than four maritime assets and two aeronautical assets due to the number of current assets belonging to the district.  For notational purposes, a response strategy in this study is represented as (\# of Maritime Assets, \# of Aeronautical Assets).  For zones of boat/helicopter events, these are assumed to be (Boats, Helicopters) whereas, for zones of cutter/airplane events, they are assumed to be (Cutters, Airplanes).  The respective frequencies of occurrence for each response strategy was modeled by empirically-constructed probability mass functions, developed for each zone.

\subsection{Location Model for Heterogeneous Assets} \label{sec:LocModel}
Using forecast demands developed from the stochastic zonal distribution model, a prescriptive approach is necessary to locate District 14's boats, cutters, airplanes, and helicopters to reduce the expected response time for future SAR events.  To formulate the mathematical programming model for this problem, we first define the following sets, parameters, and decision variables. \\\noindent 

\noindent\textbf{Sets:}
    \begin{itemize}
        \item $n \in N$: Set of asset categories, where N=\{Boat, Cutter, Airplanes, Helicopters\} for the District 14 asset location problem.
        \item $m \in M$: Set of location categories, where M=\{Harbor, Airport\} for the District 14 asset location problem.
        \item $(n,m) \in P$: Set of \textit{infeasible} asset-location combinations.  As would be expected, P=\{(Cutter, Airport), (Boat, Airport), (Airplane, Harbor), (Helicopter, Harbor)\}.
        \item $h \in H_n$: Set of all individual assets within asset category $n$.  For the problem considered herein, this set is informed by current District 14 operational capabilities.
        \item $i \in I_m$: Set of candidate homeports.  This set is informed by airports and ports identified in coordination with District 14 operational staff.
        \item $j \in J$: Set of demand nodes.  This set is informed via the historical data and the process described in Section \ref{sec:ModelParam}.
    \end{itemize}

\noindent\textbf{Parameters:}
    \begin{itemize}
        \item $c_{hi}$: Time (hours) to (re)assign asset $h \in H_n$ to candidate homeport $i \in I_m$.  Each $c_{hi}$-value is computed by dividing the distance (nmi) between the current homeport of asset $h\in H_n$ and candidate homeport $i\in I_m$ by the cruise speed (knots) of asset $h\in H_n$.  Note that assigning an asset to its current homeport yields $c_{hi}=0$.
        \item $d_{hij}$: Time (hours) to deploy asset $h \in H_n$ from candidate homeport $i \in I_m$ to demand node $j \in J$.  The time to deploy assets is computed by dividing the distance (nmi) between the candidate homeport $i\in I_m$ and the demand node $j\in J$ by the maximum speed (knots) of asset $h\in H_n$.
        \item $l_{nj}$: Level of demand (number of assets required) for asset type $n \in N$ at site $j \in J$. 
        \item $u_h$: Monthly hours allocated for SAR operations by asset $h \in H_n$. This allocation is user-determined and depends on organizational staffing, training, maintenance, and budget limitations.
        \item $t$: Time (hours) required on-site to complete a SAR mission.  We assume a constant value of $t=1.5$ for all SAR missions, but a practitioner can readily alter the model to accommodate alternative assumptions. 
        \item $Q$: An arbitrarily large number used within a logical constraint to relate location and allocation decisions.
    \end{itemize}

\noindent\textbf{Decision Variables:}	 
    \begin{itemize}
        \item $x_{hi}$: Binary variable equal to 1 if asset $h \in H_n$ is assigned to candidate homeport $i \in I_m$, and 0 otherwise.
        \item $y_{hij}$: Integer-valued number of SAR events to which asset $h \in H_n$, assigned to candidate homeport $i \in I_m$, is poised to respond at demand nodes $j \in J$.
    \end{itemize}

Given this framework, we propose the following formulation for the Optimal Location Model for Heterogeneous Search and Rescue Assets:

{\allowdisplaybreaks
\begin{align}
     \label{Con:ObjFun1} \min_{\boldsymbol{x}, \boldsymbol{y}} & \left(f_1(\boldsymbol{x}), f_2(\boldsymbol{y}) \right) \\
     \label{Con:WeightedSum1} \text{s.t. }	& f_1(\boldsymbol{x})=\sum_{(m,n)\in \{M\times N \}\setminus P} \left( \sum_{h \in H_n} \sum_{i \in I_m} c_{hi}x_{hi}\right),\\
     \label{Con:WeightedSum2}  &  f_2(\boldsymbol{y})=\sum_{(n,m)\in \{N\times M \}\setminus P} \left( \sum_{h \in H_n} \sum_{i \in I_m} \sum_{j \in J} d_{hij}y_{hij}\right), \\
     \label{Con:AssignOne} & \sum_{m\in M} \sum_{i \in I_m} x_{hi} = 1, \ \forall \ h \in H_n, n \in N,\\
     \label{Con:BlockMismatch} & \sum_{h \in H_n} \sum_{i \in I_m} x_{hi} = 0,\ \forall \ (m,n) \in P, \\
     \label{Con:MeetDemand} & \sum_{h\in H_n} \sum_{m\in M} \sum_{i\in I_m} y_{hij} \geq l_{nj}, \ \forall \ n\in N, j\in J,\\
     \label{Con:PoiseAssets} & \sum_{j\in J} y_{hij} \leq Q x_{hi}, \ \forall \ h \in H_n, n \in N, i\in I_m, m\in M,\\
     \label{Con:OpCapacity} & \sum_{m\in M} \sum_{i \in I_m} \sum_{j \in J} \left(2d_{hij} + t \right)y_{hij} \leq u_h,\ \forall \ h \in H_n, n \in N, \\
      \label{Con:Binary} & x_{hi} \in \{0, 1\},\ \forall \ h \in H_n, n \in N, i \in I_m, m \in M,\\
      \label{Con:Integer} & y_{hij} \in \mathbb{Z}_+,\ \forall \ h \in H_n, n \in N, i \in I_m, j \in J, m \in M.
\end{align}}

The multiobjective formulation seeks to minimize two objectives.  The first objective \eqref{Con:WeightedSum1} calculates the total time (in hours) required to reassign assets to new homeports, and the second objective \eqref{Con:WeightedSum2} computes the expected total time (in hours) to deploy assets from their homeports to assigned demand nodes.  Although the two objectives share the same units of measure, this is merely a coincidence resulting from the decision to model asset relocation costs in units of time to physically move the assets between locations.  Such an objective is a proxy for the costs of moving assets, a measure that is more elusive to compute due to the imposition of both tangible costs (e.g., purchasing real estate, leasing warehousing for aircraft or berthing for seagoing vessels, building support infrastructure) and intangible costs (e.g., political negotiation, military appropriations testimony) when relocating assets.  Moreover, the two objectives are not necessarily well-scaled, although we find them to be well-scaled for the District 14 problem instances during testing in Section \ref{sec:TRA}.

Constraint \eqref{Con:AssignOne} ensures that every asset is assigned to exactly one homeport, whereas Constraint \eqref{Con:BlockMismatch} prevents erroneous assignments, such as homeporting a maritime asset at an airport.  Constraint \eqref{Con:MeetDemand} requires the demand requirement for each asset type to be met.  Constraint \eqref{Con:PoiseAssets} only allows assets to be deployed from a homeport if they are assigned to that homeport.  An arbitrarily large value is used for Q (e.g., $Q=100 |J|$), indicating that assigned assets can respond to at most $100 |J|$ SAR events.  Constraint \eqref{Con:OpCapacity} enforces the monthly utilization rates for individual assets, where each SAR event imposes the time to travel to and from the demand node as well as the duration of the mission.  Constraints \eqref{Con:Binary} and \eqref{Con:Integer} respectively enforce binary and non-negative integer constraints on the location- and assignment-related decision variables.

Among the methods available to address the multi-objective nature of the formulation \cite{ehrgott2005multicriteria}, we utilize the Weighted Sum Method for testing.  By comparison among the more commonly used techniques, the objective function metrics are not particularly well suited for the $\varepsilon$-constraint Method, and we seek to explore the Pareto optimal set of solutions to a greater degree than would be afforded by the Lexicographic Method.  Hereafter, we denote $w_1$ and $w_2$ to be the respective non-negative weights on objectives $f_1(\boldsymbol{x})$ and $f_2(\boldsymbol{y})$ when utilizing the Weighted Sum Method, where $w_1+w_2=1$.

This formulation can be used to consider a number of scenarios.  Within this study, we consider how District 14 should posture its assets such that the total response is minimized in the following scenarios:
\begin{enumerate}
    \item District 14 is restricted to operating only out of its currently owned stations.
    \item District 14 is able to position its assets at any operational airport or harbor across the AOR.
    \item Each of Scenarios 1 and 2, when District 14 is simultaneously seeking to reduce the cost associated with reassigning assets to new homeports.
\end{enumerate}

Alterations can be made to the sets and parameters to consider each of these scenarios.  The homeports $i \in I_m$ are adjusted to model, alternatively, only locations currently operated by District 14 or a large collection of airports and harbors across the region.  Similarly, the objective function weights $(w_1,w_2)$ are adjusted based on the relative priority for the second objective function for either of the Scenario 3 variants.

To identify an optimal solution to this integer linear programming formulation, there exist a plethora of readily available, commercial solvers (e.g, CbC, CPLEX, FICO-Xpress, gurobi, MOSEK).  To encode the formulation, many algebraic modeling languages (e.g., AMPL, GAMS, LINGO, OPL) and general programming languages (e.g., C++, Java, Python, R) likewise abound.  However, we also required a modeling environment within which we could automate the computation of several complex parameters to calculate instance parameters.  Among the various parametric calculations, we needed to quickly compute the haversine distances between all considered locations, and the corresponding travel times for each asset based on vehicle specifications.  Based on the simplicity of coding environment and readily available trigonometric functions, we modeled the formulation using General Algebraic Modeling Software (GAMS) version 25.1.3 and invoked CPLEX version 12.7 to identify solutions to specific instances. 

\section{Testing, Results, and Analysis} \label{sec:TRA}
Within this section, we provide in Section \ref{sec:Data} an overview of the historical SAR event data set and data cleaning requirements, present in Section \ref{sec:SZDMresults} the results of the stochastic zonal distribution model, and report and discuss in Section \ref{sec:LocResultsAndAnalysis} the results of solving the prescriptive location model for District 14.  In subsequent examinations in the final subsection, we also apply the Weighted Sum Method to identify non-dominated (i.e., Pareto optimal) solutions for both variants of Scenario 3 under both the 50th and 75th percentile realizations of expected demand levels, followed by an examination of the relative sensitivity of the enterprise performance to this assumption.

\subsection{Historical Demand Data} \label{sec:Data}
For this study, the historical SAR data was provided from the USCG Marine Information for Safety and Law Enforcement (MISLE) database by the USCG Research and Development Center.  The event data consisted of a list of SAR events spanning from December 2010 through May 2018.
    
\begin{figure}[!ht]
\centering
\includegraphics[width=0.85\textwidth]{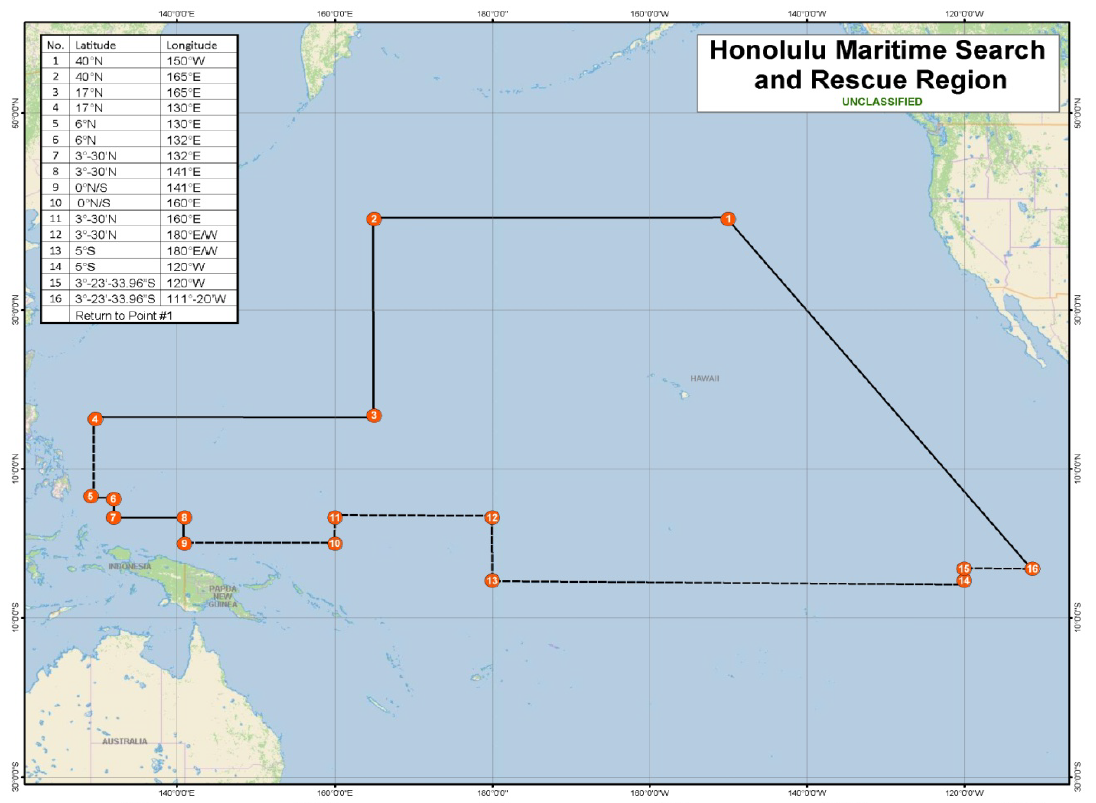}
\caption{Honolulu Maritime SAR Region \cite{CG14SARPlan2014}}
\label{fig:SAR_Region}
\end{figure}

Not every event within the data was applicable to this study, for which the scope is limited to SAR events within District 14's AOR that require the response of USCG assets.  Figure \ref{fig:SAR_Region} depicts the Honolulu Maritime Search and Rescue Region, as stipulated in the district's search and rescue plan \cite{CG14SARPlan2014} and used as a SAR-specific proxy for the (larger) District 14 AOR for this analysis.  Thus, this analysis supports asset location planning for response to events within the SAR region supported by deliberate district plans.

Three categorical subsets of data were removed from the data set.  First, SAR event records having the subtype `MEDICO' were removed because they consist of medical consultation only via phone or radio communication and do not require the response of USCG assets.  The second and third categories of event records removed from the data set were SAR events for which, respectively, either no GPS location was provided with the SAR event record or the event occurred outside the region depicted in Figure \ref{fig:SAR_Region}.   

There did exist SAR event records within the MISLE data that occurred within the Honolulu Maritime SAR Region but for which no assets were dispatched in response.  These records did not exhibit any geographic or subtype commonality, and it is assumed the records were correct and the decision not to dispatch assets in response to these events was deliberate (e.g., support was provided by an non-USCG asset).  As such, these SAR event records were retained within the data set.  Table \ref{tab:DataClean} summarizes the results of the data cleaning; 91.52\% of the historical record data was retained to inform the stochastic zonal distribution model.

\begin{table}[!ht]
    \centering
    \small
    \caption{Data Cleaning Summary for SAR event Records, Dec 2010 - May 2018}
    \label{tab:DataClean}
        \begin{tabular}{lr} 
         \textbf{Category} & \textbf{No. Records} \\  \hline
         Initial Data Set & 4315  \\ 
         MEDICO events & 90  \\
         Missing GPS Data & 38  \\
         Outside SAR Region & 238  \\ 
         Final Data Set & 3949  \\ 
        \end{tabular}
\end{table}

Providing response to forecast events are District 14's fleet of 21 assets: eight boats, nine cutters, two helicopters, and two airplanes.  The district operates out of six locations: four harbors and two air stations.  Expanding consideration for possible asset locations across the Pacific region yields 50 civilian locations (i.e., 28 harbors and 22 airports) identified for this study.

Several assumptions were necessary to leverage the MISLE data in both the forecasting and optimization models.  First, the historical SAR data is assumed to be indicative of future trends in USCG District 14's SAR region.  Second, the 3949 records in the cleaned data set are assumed to accurately represent the SAR missions within the District 14 SAR region and the corresponding demands upon District 14 assets from December 2010 through May 2018.  The last assumption relates to the distance measurements used in the location problem formulation.  This study uses the haversine formula to calculate distances between candidate homeports and both other homeports and demand nodes.  Using this distance calculation assumes both that the Earth is perfectly spherical, an acceptable approximation, and USCG assets travel by the distance-minimizing route to and from SAR events.  Thus, the calculations neglect the additional distances to route assets around islands as well as any deviations from the most direct route due to asset traffic, tides, or weather.

\subsection{Stochastic Zonal Distribution Model Results} \label{sec:SZDMresults}
As introduced in Section \ref{sec:ModelParam}, we use the stochastic zonal distribution model to sequentially construct SAR demand zones, characterize the stochastic nature of SAR event demands within each zone, and identify the response levels of such SAR demands, in terms of the type and numbers of the various heterogeneous assets.

\subsection*{Step 1.  Hierarchical Clustering of Historical SAR Data}
District 14 is comprised of three response organizations among the MISLE data:  Sector Guam, Sector Honolulu, and District 14 Headquarters.  As such, all SAR events were first sorted among four categories: (1) Sector Guam boat/helicopter events, (2) Sector Guam cutter/airplane events, (3) Sector Honolulu/District 14 Headquarters boat/helicopter events, and (4) Sector Honolulu/District 14 Headquarters cutter/airplane events.  Sector Honolulu and District 14 Headquarters were combined because of a significant level of overlap between the geographic dispositions of historical SAR events to which they respectively responded.  

The locations of the SAR events were represented on a two-dimensional Cartesian plane using their respective longitude and latitude coordinates, adjusting coordinates for the presence of the antimeridian to reference all coordinates via a common direction (i.e., west) with respect to the prime meridian.  Although such a projection is not equivalent to the geometry along the surface of the Earth, we assume the error to be negligible, i.e., we assume SAR events clustered together on a planar projection would be clustered likewise on the surface of the a sphere.  Alternatively, the same clustering procedures were implemented on the data when represented in a three-dimensional Cartesian space (i.e., considering direct, under-the-ocean-surface Euclidean distances), attaining similar results.

To determine the number of clusters within each of the four categories, a series of elbow curves were generated in Python.  A noted limitation of the $k$-Means approach is the sensitivity of the procedure to the selection of starting centroid locations; ill-selected starting locations can lead to poor clustering results.  A solution proposed by Arthur and Vassilvitskii \cite{arthur2007k} is the $k$-Means++ approach, wherein the initial centroid point is randomly selected, and subsequent centroid points are selected based on a probability that is a function of the shortest distance between the proposed center point and previously selected center points.  Once the $k$ initial center points are selected, the traditional $k$-Means procedure is implemented to cluster the data points.  The $k$-Means++ method of selecting initial center points was utilized for the clustering in this study.
    
\begin{figure}[!ht]
\centering
\includegraphics[width=0.85\textwidth]{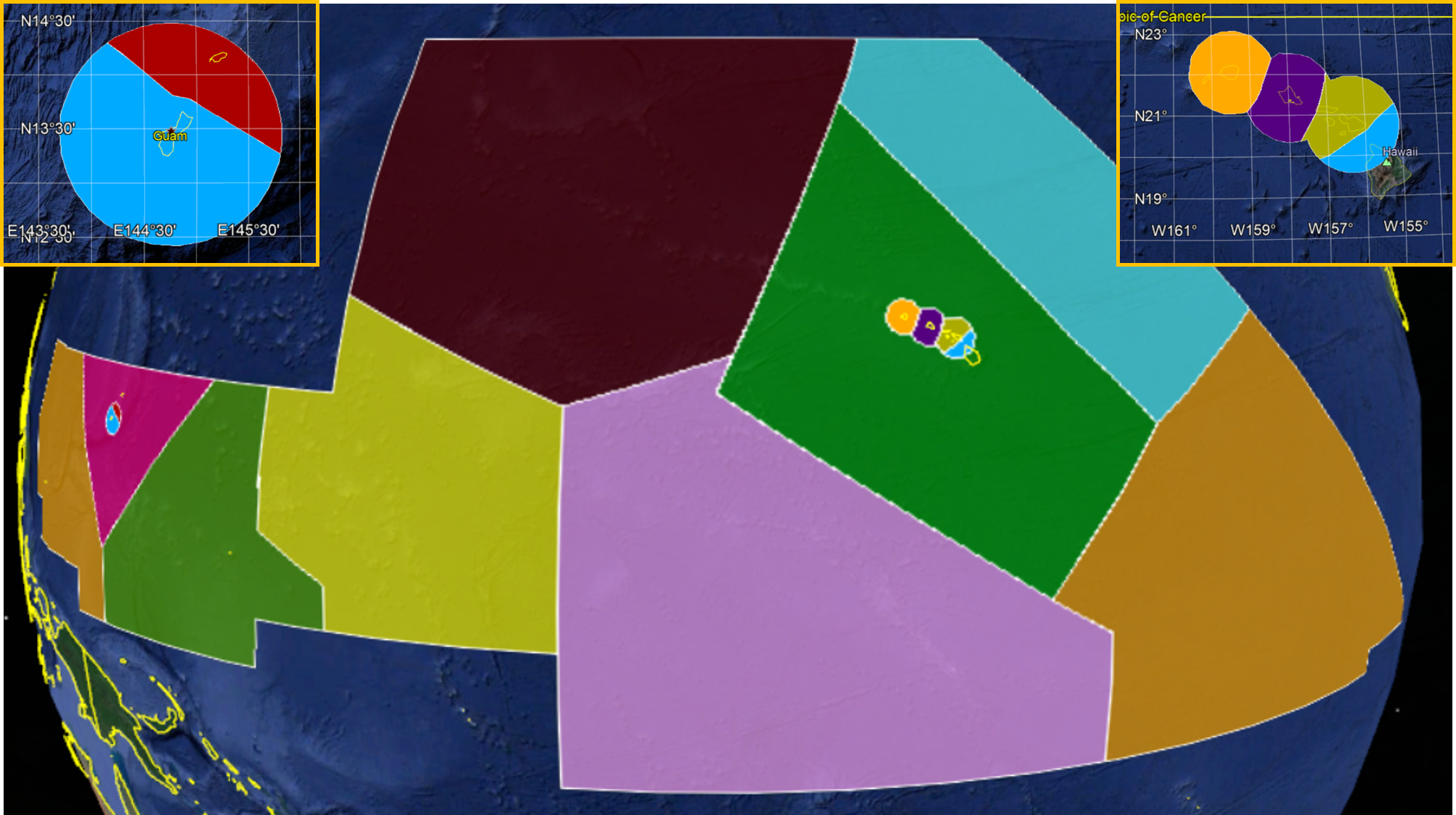}
\caption{Geographical Approximation of the 15 SAR event zones for the District 14 Stochastic Zonal Distribution Model}
\label{fig:Zones}
\end{figure}

Executing this clustering procedure, 15 clusters are generated.  Of these 15 clusters, six clusters consisted of boat/helicopter events (two surrounding Guam and four surrounding the Hawaiian Islands), and nine clusters consisted of cutter/airplane events (three within Sector Guam's AOR and six within Sector Honolulu/District 14 Headquarter's AOR).  Zones boundaries around the clusters were identified to partition the District 14 AOR, as depicted in Figure \ref{fig:Zones}.

The zonal distribution model established by Azofra et al. \cite{azofra2007optimum} utilizes `superaccident' sites to represent the aggregated demand node of SAR events within each zone.  The authors' model computes the location of these superaccidents via the arithmetic mean of the longitudes and latitudes for all events within the respective zones. Razi and Karatas \cite{razi2016multi} also leveraged the concept of superaccident sites, but improved upon the work of Azofra et al. \cite{azofra2007optimum} by accounting for varying weights of events; they implemented a weighted $k$-Means clustering algorithm to account for different event magnitudes in the zonal grouping of SAR events.  We adopt a similar approach and calculate the superaccidents for all 15 clusters using weighted event data, establishing the magnitude of each SAR event as the total number of activities associated with the event's case file.  An activity is created for each resource sortie assigned to the case or whenever the nature of a case has changed.  The total number of activities was concluded to be an appropriate measure of SAR event magnitude, based upon the assumption that SAR events which are larger in scale are inclined to require more response assets.  The weighted centroids of each cluster were computed, and these locations were designated as superaccident sites.  Table \ref{tab:Weight Comp} presents the superaccident sites identified in this study.

\begin{table}[!ht]
    \centering
    \small
    \caption{Weighted Superaccident Locations for the District 14 Stochastic Zonal Distribution Model}
    \label{tab:Weight Comp}
    \begin{tabular}{l|l}
       \textbf{Zone} & \textbf{Coordinates} \\ \hline
       Guam-0 & $13^{\circ}25'52.179''$N, $144^{\circ}41'45.1104''$E \\
       Guam-1 & $14^{\circ}2'25.3788''$N, $145^{\circ}9'29.6346''$E \\
       Hawaii-2 & $20^{\circ}48'36.1902''$N, $156^{\circ}35'53.9484''$W \\
       Hawaii-3 & $21^{\circ}59'35.8044''$N, $159^{\circ}24'40.143''$W \\
       Hawaii-4 & $21^{\circ}22'46.7472''$N, $157^{\circ}56'10.5138''$W \\
       Hawaii-5 & $19^{\circ}58'4.332''$N, $156^{\circ}0'57.852''$W \\
       Guam-6 & $14^{\circ}27'44.985''$N, $145^{\circ}36'1.0074''$E \\
       Guam-7 & $8^{\circ}48'38.415''$N, $136^{\circ}29'7.35''$E \\
       Guam-8 & $6^{\circ}59'52.548''$N, $153^{\circ}28'20.3874''$E \\
       Hawaii-9 & $20^{\circ}3'10.011''$N, $156^{\circ}35'43.8144''$W \\
       Hawaii-10 & $29^{\circ}42'5.022''$N, $179^{\circ}43'38.6322''$E \\
       Hawaii-11 & $8^{\circ}4'46.416''$N, $131^{\circ}15'44.3736''$W \\
       Hawaii-12 & $10^{\circ}41'58.4484''$N, $168^{\circ}42'48.7764''$E \\
       Hawaii-13 & $3^{\circ}46'51.2034''$N, $164^{\circ}55'1.362''$W \\
       Hawaii-14 & $27^{\circ}48'16.1058''$N, $147^{\circ}54'4.2366''$W \\
    \end{tabular}
\end{table}

\subsection*{Step 2.  Probabilistic Representation of SAR Demand for each Zone}
As discussed in Section \ref{sec:ModelParam}, previous studies of SAR operations consider these events to be Poisson distributed.  Reviewing the requirements for a Poisson process (e.g., see \cite{ross2014introduction}), historic SAR events largely adhere to these criteria.  

Of note, the adoption of a Poisson distribution requires the events to be both independent and stationary.  Independence suggests that the probability of one event occurring does not effect the probability of another event occurring.  We assume SAR events to be independent, in general, although we recognize that there do exist instances when this may not hold.

Stationarity requires that, for any interval in the considered time series, the rate of arrival remains constant.  A frequent challenge to this assumption in Coast Guard studies results from seasonal trends in the event data.  The determination of whether the SAR data exhibited statistically significant levels of seasonality was made by examining autocorrelation function plots generated using the time series analysis functions within JMP.  The autocorrelation for a lag-time $k$ is given by Equation \ref{eq:Autocorrelation}, where $y_t$ represents the number of SAR events at time $t$.  An autocorrelation of $|r_k| \approx 1$ is representative of a larger relationship, whereas an autocorrelation near 0 is indicative of little-to-no relationship between the data.  In particular, the existence of recurring annual levels of seasonality would result in larger levels of autocorrelation near lag-times in increments of 12.  
\begin{equation} \label{eq:Autocorrelation} r_k = \frac{\displaystyle\sum^N_{t = k + 1}\left(y_t - \bar{y}\right)\left(y_{t - k} - \bar{y}\right)}{\left(y_t - \bar{y}\right)^2} \end{equation}

Although fluctuations in levels do occur, a visual review of the autocorrelation plots along with the autocorrelations for lag-times of 12- and 24-months, as reported in Table \ref{tab:AutoCor}, provide no statistically significant evidence of a relationship between the levels of SAR activity.  This lack of exhibited autocorrelation was further examined by a visual inspection of time series plots.  The time series plots suggested that SAR events in the Pacific Ocean region do not exhibit seasonality, with the exception of a season spike in workload in the zone encompassing O'ahu during the summer months, particularly in July.

\begin{table}[!ht]
    \centering
    \small
    \caption{Autocorrelations for each SAR event Cluster at 12-month and 24-month Lag Times}
    \label{tab:AutoCor}
    \begin{tabular}{l|rr}
       \textbf{Zone}  & \textbf{12 Month} & \textbf{24 Month} \\  \hline
       Guam-0 & 0.0284 & -0.0506 \\
       Guam-1 & 0.0873 & -0.0679 \\
       Hawaii-2 & 0.0957 & -0.0561 \\
       Hawaii-3 & -0.2557 & -0.0172 \\
       Hawaii-4 & 0.1302 & 0.0377 \\
       Hawaii-5 & -0.1089 & -0.2341 \\
       Guam-6 & 0.0092 & 0.1627 \\
       Guam-7 & -0.0548 & -0.1102 \\
       Guam-8 & -0.1302 & 0.0589 \\
       Hawaii-9 & 0.1122 & -0.0357 \\
       Hawaii-10 & -0.0332 & -0.0327 \\
       Hawaii-11 & 0.0331 & -0.0074 \\
       Hawaii-12 & -0.0232 & -0.1070 \\
       Hawaii-13 & -0.0219 & -0.0462 \\
       Hawaii-14 & -0.0519 & 0.2372 \\  \hline
       All Clusters & 0.0808 & -0.0940 \\
    \end{tabular}
\end{table}

Data is considered to be stationary if, for any time interval in the series, the distribution remains the same.  From Table \ref{tab:AutoCor}, we note that all the autocorrelation levels are fairly close to zero, with only the following clusters exhibiting $|r_k| > 0.1$ for 12-month lag times: Hawaii-3, Hawaii-4, Hawaii-5, Guam-8, and Hawaii-9.  We conclude that the preponderance of the zonal SAR data is stationary, though there may be noticeable fluctuation in the rate of events for the excepted zones.  Thus, we further conclude that not all zones were shown to be strictly Poisson-distributed because of these non-stationary elements.
 
To circumvent the non-stationary aspects of these clusters, Gamma-Poisson distributions were fit to the monthly time-series data.  Gamma-Poisson is a mixture of the two probability distributions, specifically for models where the count of individual events $x$ is Poisson distributed with a rate $\lambda$, which is itself Gamma distributed with parameters $\alpha$ and $\beta$.  In other words, the count of events $x$ is conditional on the Poisson parameter $\lambda$, and $\lambda$ is conditional on the Gamma parameters $\alpha$ and $\beta$.
    
Utilizing the \textit{Distribution} application in JMP, both Poisson and Gamma-Poisson distributions were constructed for the monthly time-series data in each of the zones, and a Pearson's chi-squared test was applied to evaluate the goodness-of-fit.  The parameter estimates and corresponding goodness-of-fit \textit{p}-values are displayed in Table \ref{tab:DistParam}.  It was found that modeling the data with the Gamma-Poisson distribution resulted in better fits for all zones except for Hawaii-2 and Hawaii-13.  In these zones, there wasn't evidence of a sufficient variation in $\lambda$ to merit the use of a Gamma-Poisson distribution; these clusters were henceforth modeled by a Poisson distribution.

\begin{table}[!ht]
    \centering
    \small
    \caption{Summary of Pearson's Chi-Squared Goodness-of-Fit Test Results}
    \label{tab:DistParam}
    \begin{tabular}{ l |rr | rrr}
       & \multicolumn{2}{c|}{\textbf{Poisson Dist.}} & \multicolumn{3}{c}{\textbf{Gamma-Poisson Dist.}} \\
       \textbf{Zone}  &$p$-Value& $\lambda$ & $p$-Value & $\alpha$ & $\beta$ \\ \hline
       Guam-0 & 0.2191 & 5.433 & 0.4575 & 52.748 & 0.103 \\
       Guam-1 & 0.1635 & 0.533 & 0.4458 & 4.069 & 0.131 \\
       Hawaii-2 & 0.8786 & 6.256 & -- & -- & -- \\
       Hawaii-3 & 0.0126 & 3.044 & 0.4535 & 8.672 & 0.351 \\
       Hawaii-4 & 0.0101 & 13.022 & 0.3876 & 39.105 & 0.333 \\
       Hawaii-5 & 0.1137 & 1.378 & 0.4449 & 8.202 & 0.168 \\
       Guam-6 & 0.0605 & 2.689 & 0.4422 & 11.951 & 0.225 \\
       Guam-7 & 0.3751 & 1.756 & 0.4624 & 51.647 & 0.034 \\
       Guam-8 & 0.0325 & 1.733 & 0.3485 & 7.535 & 0.230 \\
       Hawaii-9 & 0.2693 & 3.644 & 0.4438 & 50.611 & 0.072 \\
       Hawaii-10 & 0.0053 & 0.978 & 0.4006 & 2.560 & 0.382 \\
       Hawaii-11 & 0.0775 & 0.833 & 0.5569 & 3.254 & 0.256 \\
       Hawaii-12 & 0.0307 & 0.700 & 0.5664 & 2.047 & 0.342 \\
       Hawaii-13 & 0.6819 & 0.833 & -- & -- & -- \\
       Hawaii-14 & 0.0450 & 1.044 & 0.4248 & 4.332 & 0.241 \\
    \end{tabular}
\end{table}

\subsection*{Step 3.  Probabilistic Representation of SAR event Asset Response Requirements}
In the MISLE data, there was an accompanying record of the assets dispatched to each SAR event.  Among the data set were 1133 uniquely named assets that supported SAR responses over the approximately seven-year period.  Some of these assets belonged to District 14, but most were other military, commercial, and private boats and aircraft.  From an modeling perspective, the presence of these non-USCG assets cannot be guaranteed when a SAR event occurs, yet District 14 still has the responsibility to coordinate the SAR responses.  Adopting a conservative view, the interpretation of these asset records focused less on the specific assets that participated in the SAR response and more on the general type and number of assets required to respond.
    
To gain insight into these general strategies of SAR operations, the 1133 uniquely named assets were categorized as either aeronautical or maritime assets.  From the synthesis of these event records, three general strategies regarding SAR operations emerged: respond with  aeronautical assets only, respond with maritime assets only, or respond with a combination of maritime and  aeronautical assets.  A fourth type of response in which no assets were dispatched was also observed. These events remained in the data set to inform the forecast for future demand locations but not the level of response required for SAR events, lest we accidentally assume away potential demands for District 14 asset response.   Table \ref{tab:SARStrat2} presents the relative response strategy rates for each zone.

\begin{table}[!ht]
    \centering
    \small
    \caption{Response Strategy Proportions (\%), by Zone, for the District 14 Stochastic Zonal Distribution Model}
    \label{tab:SARStrat2}
    \begin{tabular}{ l | rrrr }
        & \textbf{Aircraft} & \textbf{Maritime} & \textbf{Maritime \&} \\
       \textbf{Zone} & \textbf{Only} & \textbf{Only} & \textbf{Aircraft} \\ \hline
       Guam-0 & 8.642 & 73.827  & 17.531  \\
       Guam-1 & 59.524  & 33.333  & 7.143  \\
       Hawaii-2 & 13.586  & 67.261  & 19.154  \\
       Hawaii-3 & 18.386  & 53.812  & 27.803  \\
       Hawaii-4 & 28.066  & 50.118  & 21.816  \\
       Hawaii-5 & 44.928  & 13.043  & 42.029  \\
       Guam-6 & 23.636  & 56.364  & 20.000  \\
       Guam-7 & 1.220  & 91.463  & 7.317  \\
       Guam-8 & 11.111  & 61.111  & 27.778 \\
       Hawaii-9 & 48.691  & 22.513  & 28.796  \\
       Hawaii-10 & 32.000  & 54.000  & 14.000  \\
       Hawaii-11 & 9.375  & 75.000  & 15.625  \\
       Hawaii-12 & 28.947  & 36.842  & 34.211  \\
       Hawaii-13 & 40.000  & 37.778  & 22.222  \\
       Hawaii-14 & 30.357  & 50.000  & 19.643  \\
    \end{tabular}
\end{table}  
    
Having determined the tendencies for District 14 to respond to various SAR events with these three general strategies, it was necessary to identify frequency of the different response volumes, by asset type, to SAR events for each strategy.  Since the deterministic location models for this thesis only considered the allocation of USCG assets, the number of  aeronautical assets considered was limited to two and the number of maritime assets considered was limited to four.  If, for instance, a historical SAR event was supported by six maritime assets, that would be considered in this analysis as a four maritime asset case.  To complete the third step of the stochastic zonal distribution model, we computed each of the conditional probabilities of each level of response for events, by zone, which we forgo presenting herein for the sake of brevity.  

As an illustrative example, consider zone Guam-0.  The stochastic zonal distribution model identified Guam-0 as a cluster of boat-supported SAR events coordinated by Sector Guam.  The center of SAR operations for the zone was found to be $13^{\circ}25'52.179''$N, $144^{\circ}41'45.1104''$E, as reported in Table \ref{tab:Weight Comp}.  The emergence of SAR events within the zone is modeled with the Gamma-Poisson probability distribution having $(\lambda,\alpha,\beta)=(5.433,52.748, 0.103)$, as presented in Table \ref{tab:DistParam}.  From Table \ref{tab:SARStrat2} and the aforementioned conditional probability calculations, SAR events in this zone will require, e.g., 0 boats and 1 helicopter for 8.4\% of the events; 0 boats and 2 helicopters for 0.2\% of the events, and 1 boat and 0 helicopters for 60.9\% of the events.

The outputs of the stochastic zonal distribution model informed a Monte Carlo simulation to replicate 10,000 months of activity for each of the 15 zones.  For each zone, 10,000 months were simulated as follows: the number of SAR events for the zone in a given month was selected based upon the corresponding distribution in Table \ref{tab:DistParam} and, for each event in that month, a SAR response was simulated using the respective probabilities.  The resulting 10,000 months of simulated data were summarized by descriptive statistics: the mean, standard deviation, and percentiles (minimum, 25\%, 50\%, 75\%, maximum).  These percentile values can serve as the relative demand volumes for each asset type in each zone, depending on the desired level of robustness.

\subsection{Location Model Results and Analysis} \label{sec:LocResultsAndAnalysis}
Having parameterized the asset demand levels for each zone, the integer linear programming formulation set forth in Section \ref{sec:LocModel} can be used to solve Scenarios 1-3 with appropriate, respective sets of objective function weights.  We solve the model for each scenario using, alternatively, the 50th and 75th percentile levels of SAR asset demands for each zone from the Monte Carlo simulation.  Whereas the 50th percentile levels are used to represent a typical month of SAR operations, the 75th percentile of demand depicts an elevated operational tempo for District 14, allowing for a risk averse perspective for strategic location of assets.  Of note, the 100th percentile levels for clusters were not considered in this analysis; such an extreme situation would likely skew the results beyond what should be practically implemented.

For both the 50th and 75th percentiles of SAR demand, the instances are iteratively solved via the Weighted Sum Method, considering 0.2 increments in each of the respective objective function weights.  At the extremes, zero-valued weights were prohibited.  For example, the weights of $(w_1,w_2)=(1-\delta,\delta)$ with $\delta=0.00005$ were applied to preemptively minimize the first objective  $f_1(\boldsymbol{x})$ to maintains the current fleet posture, yet discriminate among alternative optimal solutions with respect to that objective by minimizing the second objective $f_2(\boldsymbol{y})$ to ensure that, despite being confined to their current posture, assets will respond to SAR demand via assignments to minimize the response time.

As an illustrative example that does not indicate the current locations of USCG District 14 assets (i.e., for the purpose of maintaining operational security), Table \ref{tab:ExampleSolution} presents the prescribed asset locations and event site assignments when responding to the 75th percentile demand levels, considering the larger set of possible ports and airfields at which to base assets across the Pacific Region, and when preemptively prioritizing the minimization of total expected response time (i.e., $(w_1,w_2)=(\delta,1-\delta)$.  Within Table \ref{tab:ExampleSolution}, the specific assets are denoted by a code indicating the type of asset, with a suffix of `-1', `-2', and so forth to discriminate between specific assets.  For example, 110' WPB-2 is the second of two island-class patrol boat having a length of 110 feet, and 87' CPB-1 is the first of two marine protector class patrol boats having a length of 87 feet.

\begin{table}[!ht]
    \centering
    \small
    \caption{Prescribed Asset Locations and Allocations for the 75th\% Demand Levels, Pacific Region Ports for Scenario 2}
    \label{tab:ExampleSolution}
    \begin{tabular}{c|ccc }
        \textbf{Asset Type }& \textbf{Asset} & \textbf{Homeport} & \textbf{Assigned Zones} \\ \hline
       \multirow{9}{*}{Cutter} &  225' WLB-1 & Honolulu Harbor &  --  \\
        &  225' WLB-2 &  Port of Kwajalein & Hawaii-12  \\
        &  110' WPB-1 & Tomil Harbor  & Guam-7  \\
        &  110' WPB-2 & Port of Johnson Atoll  & Hawaii-13  \\
        &  FRC-1 & Hilo Harbor  &  Hawaii-11, Hawaii-14 \\
        &  FRC-2 & Port of Kailua Kona  & Hawaii-9, Hawaii-11  \\
        &  FRC-3 & Pohnpei Harbor  &  Guam-8 \\
        &  87' CPB-1 & Port of Tinian  &  Guam-6 \\
        &  87' CPB-2 &  Port of Midway Islands &  Hawaii-10  \\ \hline
       \multirow{8}{*}{Boat} &  45' RBM-1 & Honolulu Harbor  & --  \\
        &  45' RBM-2 & Maalaea Harbor  &  Hawaii-2 \\
        &  45' RBM-3 & Nawilili Harbor  &  Hawaii-3 \\
        &  45' RBM-4 & Nawilili Harbor  &  --  \\
        &  45' RBM-5 & Apra Harbor  &  Guam-0 \\
        &  45' RBM-6 & Apra Harbor  &  -- \\
        &  45' RBM-7 & Pearl Harbor   &  Hawaii-4 \\
        &  45' RBM-8 & Kawaihae Harbor  &  Hawaii-5 \\ \hline        
       \multirow{2}{*}{Airplane} &  C-130J-1 & Kona Intl Airport  & Hawaii-9, 10, 13, \& 14  \\
        & C-130J-2 & Chuuk Intl Airport  & Guam-6, Guam-8, Hawaii-12 \\ \hline
        \multirow{2}{*}{Helicopter} &  HH-65-1 & Kahului Airport  & Hawaii-2, 3, 4, \& 5  \\
        & HH-65-2 & Antonio B Won Pat Apt.  &  Guam-0, Guam-1
    \end{tabular}
\end{table}  

The solution identified in Table \ref{tab:ExampleSolution} does minimize the total expected response time to SAR events, but it does not use assets efficiently, and the prescribed solution should be tailored for implementation.  For example, two of the 45' Reponse Boat Medium (RBM) vessels are stationed in Nawilili Harbor, yet only one is assigned to respond to SAR events in the Hawaii-3 zone.  The workload would surely be shared during operations.  In general, the lack of assigning zones to 225' WLB-1, 45' RBM-1,  45' RBM-4, and 45' RBM-6 indicates, for this instance, that District 14 has more vessels than needed \textit{to attain this total expected response time} and for these relative objective function weights.  Such a conclusion does not necessarily hold under differing conditions, so we caution a reader against making generalizations.  Also neglected from this solution are the potential efficiencies in training personnel and maintaining equipment gained from co-locating similar assets; thus, an extension to this research might consider either a reduced set of locations or provide an incentive to co-locate similar assets.

Table \ref{tab:ScenarioResults} presents the optimal objective function values over the range of weights considered, considering either only current stations operated by District 14 or using candidate homeports across the Pacific region, and alternatively under either 50th or 75th percentile levels of asset demand.  The results in the first row (i.e., $(w_1,w_2)=(\delta,1-\delta)$) correspond to Scenario 1, in the last row (i.e., $(w_1,w_2)=(1-\delta,\delta)$) correspond to Scenario 2, and in the remaining rows correspond to alternative relative objective priorities for Scenario 3.

\begin{table}[!ht]
    \centering
     \caption{Optimal Objective Function Values for Alternative Weightings over Different Levels of Demand and Sets of Potential Asset Homeports}
     \small
    \label{tab:ScenarioResults}    
    \begin{tabular}{ c | r  r | r  r | r  r | r  r }
        & \multicolumn{4}{c}{\textbf{50th Percentile}} & \multicolumn{4}{c}{\textbf{75th Percentile}} \\
        & \multicolumn{2}{c|}{\textit{Current}} & \multicolumn{2}{c|}{\textit{Pacific}} & \multicolumn{2}{c|}{\textit{Current}} & \multicolumn{2}{c}{\textit{Pacific}} \\
        & \multicolumn{2}{c|}{\textit{Homeports}} & \multicolumn{2}{c|}{\textit{Region}} & \multicolumn{2}{c|}{\textit{Homeports}} & \multicolumn{2}{c}{\textit{Region}} \\        
       $(w_1, w_2)$ & $f_1(\boldsymbol{x})$ & $f_2(\boldsymbol{y})$ & $f_1(\boldsymbol{x})$ & $f_2(\boldsymbol{y})$ & $f_1(\boldsymbol{x})$ & $f_2(\boldsymbol{y})$ & $f_1(\boldsymbol{x})$ & $f_2(\boldsymbol{y})$ \\   \hline
     $(\delta,1-\delta)$ & 239.4 & 106.9 & 248.1 & 38.7 & 653.8 & 553.3 & 646.8 & 288.6 \\
        $(0.2, 0.8)$ & 2.9 & 113.9 & 116.1 & 42.9 & 141.5 & 568.1 & 197.9 & 311.9 \\
        $(0.4, 0.6)$ & 2.9 & 113.9 & 43.7 & 71.3 & 19.6 & 644.0 & 172.4 & 322.9 \\
        $(0.6, 0.4)$ & 0 & 118.2 & 4.3 & 109.5 & 2.9 & 658.5 & 115.7 & 374.5 \\
        $(0.8, 0.2)$ & 0 & 118.2 & 0.2 & 117.2 & 0 & 667.5 & 0.2 & 666.2 \\
        $(1-\delta,\delta)$ & 0 & 118.2 & 0 & 118.2 & 0 & 667.5 & 0 & 677.5 \\
    \end{tabular}
\end{table} 

From Table \ref{tab:ScenarioResults}, we note that each alternative set of weights $(w_1,w_2)$ does not yield a unique solution on the Pareto front.  For example, the set of weights for 50th percentile demand levels while considering only current homeports identifies three non-inferior solutions.  In contrast, each of the optimal solutions identified for the 75th percentile demand while considering homeports across the Pacific Region is unique; seven non-inferior solutions are identified via the Weighted Sum Method.  Should more solutions be sought, a finer resolution on weight increments is the first step towards finding them.  In general, we propose that finding a non-trivial subset of Pareto optimal solutions is practical; it enables the provision of suitable options for a decision-maker to consider rather than a single, prescribed solution.

Moreover, a comparison of the first and last rows in Table \ref{tab:ScenarioResults} illustrates the trade-off between minimizing the cost of reassigning assets and minimizing the time to respond to SAR events, i.e., the potential for improvement if District 14 is given complete latitude and support for repositioning its fleet of maritime and aeronautical assets.  These potentials for improvement increase when the size of the set of considered candidate homeports increases.  For example, in the case of 50th Percentile/Current Homeports, the total expected response time has the potential to decrease by 9.6\% whereas, in the case of 50th Percentile/Pacific Region, the total expected response time has the potential to decrease by 67.3\%.

The rationale for considering the 75th percentile in addition to the 50th percentile demand levels was to provide strategic insight for decision makers that considered the probabilistic risk of occasionally high demand months for various clusters.  This rationale leads one to question how posturing the fleet of assets with a more risk averse view will affect steady-state performance.  For both the Current Homeports and Pacific Region options for basing District 14 assets, we solved the instance for the 75th percentile demands levels with $(w_1,w_2)=(1-\delta,\delta)$), affixed the asset locations, and then resolved the instance (i.e., to determine event-to-asset assignments) to calculate the total expected response time.  We likewise affixed the optimal asset location solution for the 50th percentile demand levels and evaluated them under conditions of 75th percentile demand levels.  For both sets of basing options, the 50th and 75th demand level solutions were robust, yielding the same performance when the demand was altered as the SAR enterprise when deliberately designed for the alternative level of demand.  We thus conclude that the District 14 asset location problem is relatively insensitive to the levels of demand forecast by the Monte Carlo simulation; the priorities for asset location are consistent when total expected response time is minimized. We acknowledge that these results are instance-specific and recommend a practitioner consider such alternative demand level scenarios to characterize the robustness of a solution before implementing it.

\section{Conclusions and Recommendations}
In this study, we evaluated SAR operations across the Pacific Ocean for the USCG to inform the strategic posturing of District 14 maritime and  aeronautical assets in anticipation of emerging missions.  

We developed and demonstrated a two-stage modeling process to leverage historical event data and prescribe a location-and-demand-allocation solution for a heterogeneous fleet of assets to minimize a combination of the total expected response time to events and the time to relocate assets from their current bases.  To accomplish the spatiotemporal forecasting in the first stage model, we developed a stochastic zonal distribution model to forecast the location, frequency, and corresponding operational response requirements of future events throughout a region of interest.  To leverage different levels of demand generated via Monte Carlo simulation of the first stage results, we formulated an integer linear program to prescribe a location solution, subject to user-defined relative priorities among objectives, implemented via the Weighted Sum Method for multi-objective optimization.

Applying the two-stage model to over seven-and-a-half years of historical SAR event data for USCG District 14, we identified a set of non-inferior solutions under varying assumptions of demand levels and basing options.  Within these Pareto optimal solutions, as developed for both the 50th and 75th percentiles of forecast event and asset demand distributions, we identified asset location strategies that respectively yield a 9.6\% and 17.6\% increase in coverage over current asset basing when allowing locations among current homeports and airports, and respective 67.3\% and 57.4\% increases in coverage when considering a larger set of feasible basing locations.  Moreover, we illustrated via extended testing that, for the minimization of total expected response time to SAR events, the prescribed solutions for District 14 asset locations are insensitive to higher levels of SAR demand.

Several areas exist in which this research can be improved.  Future research can improve both the process and implementation of the two-stage model developed and demonstrated in this research.  A first extension should expand the scope of the stochastic zonal distribution model, particularly as it relates to the duration of SAR missions.  Whereas our research assumed a static notional mission length, more accurate historical records could be evaluated (e.g., for a different region that, perhaps, has more detailed records), and the location models could be modified accordingly.  By more accurately capturing the variability of mission length, the accuracy of projections for anticipated monthly asset utilization could be improved.

Additionally, extensions to the location model adopted in this study should include two additional objectives: balancing the utilization rates of assets based upon their assignment of SAR event zones, and co-locating similar assets to afford synergy in training, operations, and maintenance.

Finally, we recommend future efforts consider other USCG missions in the two-stage modeling process.  Previous studies seek to forecast the demand for and scheduling of USCG assets for other missions, and we conjecture that the application of our spatiotemporal forecasting model for other USCG missions could both provide meaningful insight regarding the total operational workload across District 14 and inform a more comprehensive asset location-allocation model to support operations.

\section*{Disclaimer} The views expressed in this article are those of the authors and do not reflect the official policy or position of the United States Air Force, the Department of Defense, or the United States Government.

\bibliography{mybibfile}

\end{document}